\documentclass[letterpaper, 10 pt, conference]{ieeeconf}
\IEEEoverridecommandlockouts
\usepackage{cite}
\usepackage{amsmath,amssymb,amsfonts}

\usepackage{mathtools}
\usepackage{graphicx}
\usepackage{mathbbol}
\usepackage{stfloats}
\usepackage{subfigure,color,balance,framed}
\usepackage{verbatim}
\usepackage{textcomp}
\usepackage{xcolor}
\usepackage{comment}
\usepackage{multirow}
\usepackage{makecell}
\usepackage{epsfig} 
\usepackage{tabularx}
\usepackage[inkscapelatex=false]{svg}

\usepackage{booktabs}
\usepackage{balance}
\usepackage{mathbbol}
\usepackage{algorithm}
\usepackage{algorithmicx}
\usepackage{algpseudocode}
\usepackage{mathrsfs}
\usepackage{threeparttable}

\newcommand{\method}[1]{\texttt{#1}}
\newcommand{\tr}{{{\mathsf T}}}

\usepackage{hyperref}
\hypersetup{colorlinks=true,      
linkcolor=blue,      
citecolor=blue,      
filecolor=blue,      
urlcolor=blue}
\newtheorem{theorem}{Theorem}
\newtheorem{lemma}{Lemma}

\newtheorem{definition}{Definition}
\newtheorem{remark}{Remark}

\usepackage[nameinlink]{cleveref}
\crefname{equation}{}{}
\crefname{theorem}{Theorem}{Theorems}
\crefname{corollary}{Corollary}{Corollaries}
\crefname{example}{Example}{Examples}
\crefname{assumption}{Assumption}{Assumptions}
\crefname{lemma}{Lemma}{Lemmas}
\crefname{proposition}{Proposition}{Propositions}
\crefname{figure}{Figure}{Figures}
\crefname{table}{Table}{Tables}
\crefname{fact}{Fact}{Facts}
\crefname{conjecture}{Conjecture}{Conjectures}
\crefname{section}{Section}{Sections}
\crefname{appendix}{Appendix}{Appendices}
\crefname{definition}{Definition}{Definitions}
\Crefname{equation}{}{}
\Crefname{theorem}{Theorem}{Theorems}
\Crefname{corollary}{Corollary}{Corollaries}
\Crefname{example}{Example}{Examples}
\Crefname{lemma}{Lemma}{Lemma}
\Crefname{proposition}{Proposition}{Proposition}
\Crefname{figure}{Figure}{Figures}
\Crefname{table}{Table}{Tables}
\Crefname{section}{Section}{Sections}
\Crefname{definition}{Definition}{Definitions}
\Crefname{appendix}{Appendix}{Appendices}

\newcommand{\ini}{\textnormal{ini}}
\newcommand{\col}{\textnormal{col}}
\newcommand{\f}{\textnormal{F}}

\newcommand{\p}{\textnormal{P}}
\newcommand{\D}{\textnormal{d}}
\newcommand{\C}{\textnormal{C}}
\newcommand{\h}{\textnormal{H}}

\newcommand{\z}{\textnormal{z}}
\newcommand{\Lo}{\textnormal{L}}

\newcommand{\R}{\textnormal{r}}

\DeclareMathOperator*{\argmin}{argmin}
\usepackage{siunitx}
\usepackage{threeparttable}
\usepackage{booktabs}
\def\BibTeX{{\rm B\kern-.05em{\sc i\kern-.025em b}\kern-.08em
    T\kern-.1667em\lower.7ex\hbox{E}\kern-.125emX}}

\title{\LARGE \bf Dictionary-free Koopman Predictive Control for Autonomous Vehicles in Mixed Traffic}

\author{Xu Shang$^{1}$, Zhaojian Li$^{2}$, and Yang Zheng$^{1}$
\thanks{This work is supported by NSF CMMI-2320697, NSF CMMI-2320698, and NSF CAREER 2340713.}
	\thanks{$^{1}$X. Shang and Y. Zheng are with the Department of Electrical and Computer Engineering, University of California San Diego, CA 92093, USA. (x3shang@ucsd.edu; zhengy@ucsd.edu),}
    \thanks{$^{2}$ Z. Li is with the Department of Mechanical Engineering, Michigan State University, East Lansing, MI 48824, USA. (lizhaoj1@egr.msu.edu).}}%
\begin{document}
\maketitle

\begin{abstract}
Koopman Model Predictive Control (\method{KMPC}) and Data-EnablEd Predictive Control (\method{DeePC}) use linear models to approximate nonlinear systems and integrate them with predictive control. Both approaches have recently demonstrated promising performance in controlling Connected and Autonomous Vehicles (CAVs) in mixed traffic. However, selecting appropriate lifting functions for the Koopman operator in \method{KMPC} is challenging, while the data-driven representation from Willems’ fundamental lemma in \method{DeePC} must be updated to approximate the local linearization when the equilibrium traffic state changes. In this paper, we propose a dictionary-free Koopman model predictive control (\method{DF-KMPC}) for CAV control. In particular, we first introduce a behavioral perspective to identify the optimal dictionary-free Koopman linear model. We then utilize an iterative algorithm to compute a data-driven approximation of the dictionary-free Koopman representation. Integrating this data-driven linear representation with predictive control leads to our \method{DF-KMPC}, which eliminates the need to select lifting functions and update the traffic equilibrium state. Nonlinear traffic simulations show that \method{DF-KMPC} effectively mitigates traffic waves and improves tracking performance.  
\end{abstract}


\section{Introduction}
The transition phase of mixed traffic where human-driven vehicles (HDVs) and connected and autonomous vehicles (CAVs) coexist may last for a long time \cite{zheng2020smoothing, li2022cooperative}. It~is~widely recognized that CAVs equipped with advanced controls have great potential in mitigating traffic waves and improving traffic efficiency \cite{li2022cooperative,milanes2013cooperative,  li2017dynamical} by considering the behavior of HDVs. Due to complex human driving behaviors, the mixed traffic dynamics are nonlinear \cite{zheng2020smoothing}, and its accurate dynamics are non-trivial to obtain, which complicates the design of CAV control. Recently, data-driven control approaches for controlling CAVs in mixed traffic have attracted increasing attention \cite{wang2023deep}. These approaches approximate the nonlinear system by linear representations obtained from the collected data. These linear models are then integrated with predictive control, which leads to Koopman model predictive control (\method{KMPC}) \cite{korda2018linear} and data-enabled predictive control (\method{DeePC}) \cite{coulson2019data}.

The Koopman operator theorem is originally established for autonomous systems without control \cite{koopman1931hamiltonian}. It is extended to controlled systems in \cite{korda2018linear} by considering the infinite control sequence as an extended state. One key idea is to lift the original state of the nonlinear system to a high-dimension space via lifting functions or observables. By choosing a set of proper lifting functions, the new lifted state can propagate (approximately) linearly in the high-dimension space. The obtained linear representation can be effectively integrated with predictive control to formulate \method{KMPC} \cite{mauroy2020koopman}. Although \method{KMPC} has been applied in many fields (\emph{i.e.}, robotics \cite{haggerty2023control} and mixed traffic \cite{zhan2022data}) and shown promising performance, choosing a suitable lifting function set is non-trivial, despite recent efforts in learning-based Koopman methods~\cite{shi2022deep, han2020deep}. The improper choice of observables may lead to large modeling errors \cite{haseli2021learning}.

The \method{DeePC} \cite{coulson2019data} utilizes a data-driven linear representation from Willems' fundamental lemma \cite{willems2005note} as a predictor in predictive control.  Willems' fundamental lemma is established for linear time-invariant (LTI) systems, which utilizes a rich-enough trajectory library to construct a data-driven representation. The recent \method{DeeP-LCC} \cite{wang2023deep} uses the \method{DeePC} for the Leading Crusie Control (LCC) \cite{wang2021leading} in mixed traffic and incorporates the limits on acceleration as well as the car-following spacing as the input/output constraints. Its control performance is validated from both large-scale numerical simulations and real-world experiments~\cite{wang2022implementation}. However, \method{DeeP-LCC} needs to linearize the nonlinear traffic dynamics around an equilibrium traffic state, and then obtain a data-driven representation. Thus, the approximated model in \method{DeeP-LCC} needs to be recomputed when the traffic equilibrium changes (see \cite[Section IV]{wang2022implementation} for details). 

Built on our recent advance \cite{shang2024willems}, this paper aims~to develop a \textit{Dictionary-Free} Koopman model predictive control (\method{DF-KMPC}) for CAV control in mixed traffic. Our key idea is to construct a data-driven representation for an approximated Koopman linear model of the mixed traffic system, inspired by an extended Willems' fundamental lemma in \cite{shang2024willems}. This data-driven representation requires no lifting functions and can directly adapt to varying traffic equilibrium states. In particular, we formulate the problem of choosing the optimal Koopman linear model from a behavioral perspective (\emph{i.e.}, the trajectory of the Koopman model should be close to pre-collected trajectories from the nonlinear system). We derive effective constraints from linear system identification to refine the {dictionary-free} Koopman representation, which is solved by an iterative algorithm. The resulting data-driven representation can be viewed as an approximation of the optimal Koopman linear model over the operating region (\emph{i.e.}, the data collection area) \cite{shang2024willems}. We finally integrate~the {dictionary-free} Koopman representation with predictive control which leads to \method{DF-KMPC} for CAV control in mixed traffic. 

The remainder of this paper is structured as follows.~We review nonlinear mixed traffic dynamics and its linear models from Koopman operator and Willems' fundamental lemma in \Cref{sec:preliminary}. \Cref{sec:data-driven-Koopman} presents the construction of an approximated data-driven representation for the optimal Koopman linear embedding. \Cref{sec:results} demonstrates our numerical results of \method{DF-KMPC}. We conclude the paper in \Cref{sec:conclusion}.

\section{Preliminaries and Problem Statement}
\label{sec:preliminary}
In this section, we briefly review the dynamics of the Car-Following LCC (CF-LCC) system \cite{wang2021leading} which is a small~unit in mixed traffic. We then present its (approximated) linear models constructed from the Koopman operator theorem and Willems' fundamental lemma.
\subsection{CF-LCC system}
\begin{figure}[t]
    \centering
    \includegraphics[width=0.48\textwidth]{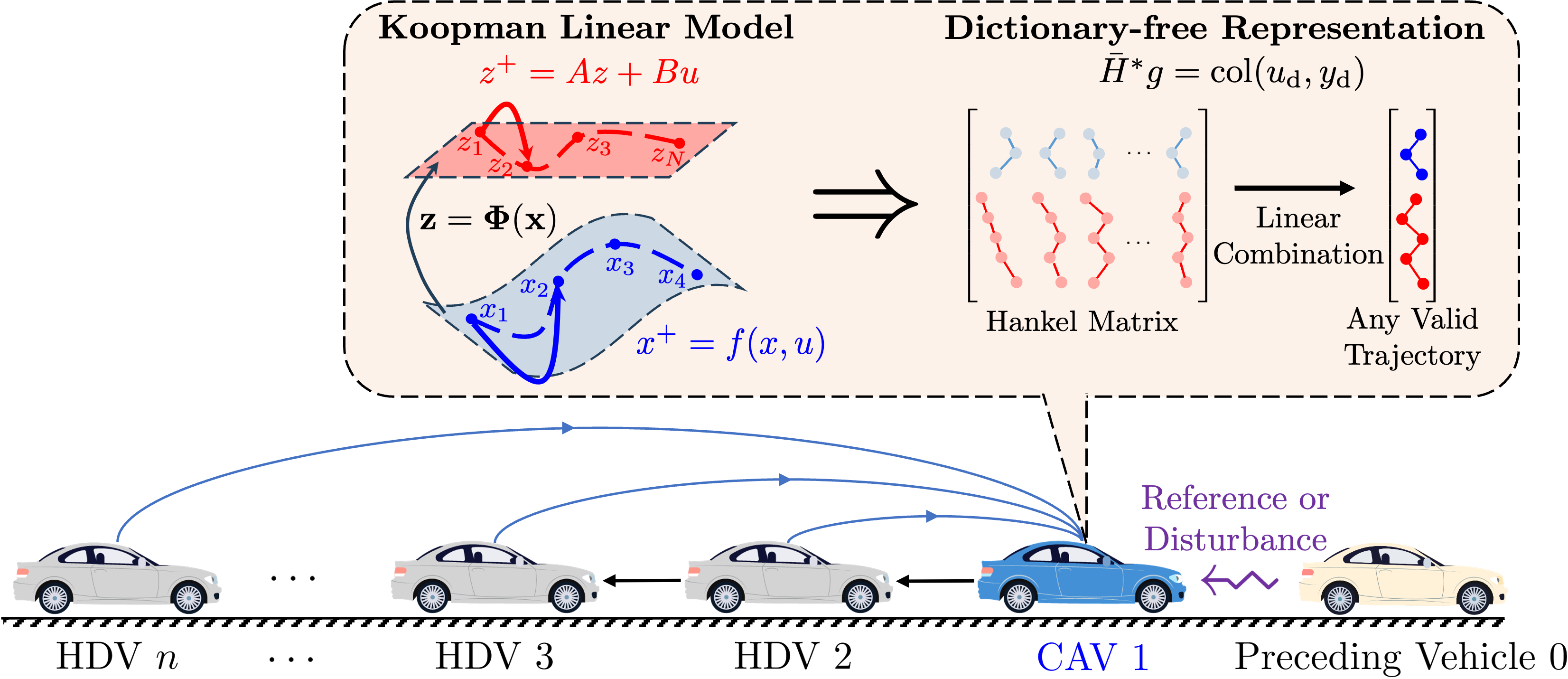}
    \caption{Schematic of a CF-LCC system. The construction of the Koopman linear model requires selection of suitable lifting functions (\emph{i.e.}, $z = \Phi(x)$), while the proposed dictionary-free representation bypasses this process.}
    \label{fig:CF-LCC}
    \vspace{-4pt}
\end{figure}
As shown in \Cref{fig:CF-LCC}, the CF-LCC considers one CAV (indexed as 0), followed by $n$ HDVs indexed as $1,\ldots, n$ from front to end. We denote the position and velocity of the $i$-th vehicle at time $k \Delta t$ where $\Delta t$ is the time interval as $p_i(k)$ and $v_i(k)$, respectively. The spacing between vehicle $i$ and its preceding vehicle is defined as $s_i(k) = p_{i-1}(k)-p_i(k)$.

We then define velocity $v_i(k)$ and spacing $s_i(k)$ as the state of the vehicle $i$ at time step $k$ and derive the state of the CF-LCC system by lumping states of all the vehicles
\[
x(k) = \begin{bmatrix}
    s_1(k), v_1(k),\ldots, s_n(k), v_n(k)
\end{bmatrix}^\tr \in \mathbb{R}^{2n}.
\]
The external input $u(k)$ includes the acceleration of the CAV and the velocity of the leading HDV, \emph{i.e.}, $u(k) = \begin{bmatrix}
    u_1(k), v_0(k)
\end{bmatrix}^\tr\! \in \! \mathbb{R}^2$.
The output of the system is identical to its state as all states are measurable. We note that most prior work on data-driven control \cite{wang2023deep, shang2024decentralized} of mixed traffic systems considers the error state, which is inconvenient as it leads to some states becoming unmeasurable and losing the spacing information of HDVs (see details in \Cref{subsec:Willems}).

Upon defining the system state, input, and output, the state-space model for the CF-LCC system can be written~as
\begin{equation}
\label{eqn:CF-LCC}
\begin{aligned}
    x(k+1) = f(x(k), u(k)), \qquad 
    y(k)  = x(k),
\end{aligned}
\end{equation}
where the nonlinear system dynamics $f:\mathbb{R}^{2n} \times \mathbb{R}^2 \rightarrow \mathbb{R}^{2n}$ is a cascading of the CAV dynamics (second-order linear model) and nonlinear HDV dynamics (\emph{e.g.}, optimal velocity model \cite{bando1995dynamical}); see more details in~\cite{wang2023deep,wang2022implementation}. An exact parametric model \eqref{eqn:CF-LCC} is non-trivial to obtain due to unknown HDV's behavior. Furthermore, even if a relatively accurate nonlinear model is identified, it may lead to a complex non-convex optimization problem when integrated with predictive control. The recently emerging data-driven methods \cite{shang2024decentralized, korda2018linear,wang2023deep} aim to derive an approximated linear model of the mixed traffic system for controller synthesis. 

\subsection{Koopman linear models for CF-LCC system}
The Koopman operator lifts the state $x(k)$ of the nonlinear system \cref{eqn:CF-LCC} to a higher-dimension space via a set of lifting functions (\emph{i.e.}, observables) \cite{korda2018linear}, where these observables propagate (approximately) linearly. Let $\phi_1(\cdot), \ldots, \phi_{n_\z}(\cdot):\mathbb{R}^{2n} \rightarrow \mathbb{R}$ be a set of linearly independent lifting functions and the output $y(k)$ is a linear map of $\Phi(x(k)):= \col(\phi_1(x(k)), \ldots, \phi_{n_z}(x(k)))$ and $u(k)$. We denote $z(k) \in \mathbb{R}^{n_\z}$ as the new lifted state $\Phi(x(k))$, which satisfies
\begin{equation}
\label{eqn:Koop-accu}
\begin{aligned}
z(k+1) & = \Phi(x(k+1)) = \Phi(f(x(k), u(k))), \\
y(k) & = Cz(k) + D u(k). 
\end{aligned}
\end{equation}
With a proper set of lifting functions, we can (approximately) represent the evolution in \cref{eqn:Koop-accu} with a parametric linear model 
\begin{equation}
    \label{eqn:Koop-linear}
    z(k+1) = Az(k) + Bu(k), \quad y(k) = Cz(k)+Du(k),
\end{equation}
where $A, B, C$ and $D$ are system matrices with  appropriate dimensions. This model \cref{eqn:Koop-linear} is known as an exact \textit{Koopman linear embedding} in \cite{shang2024willems}. Similar models as \cref{eqn:Koop-linear} are used to approximate general nonlinear dynamics in recent work~\cite{korda2018linear}. The dimension of this Koopman linear model is typically larger than the original dimension (\emph{e.g.}, $n_\z \gg 2n$ in \cref{eqn:CF-LCC}).

One can then utilize extended dynamic model decomposition (EDMD) \cite{williams2015data} to compute the data matrices $A, B, C$ and $D$ for the linear model \eqref{eqn:Koop-linear} after selecting the lifting function set $\Phi(\cdot)$. The collected input-state-output data sequence of~\cref{eqn:CF-LCC} can be arranged as 
\[
\begin{aligned}
X = \begin{bmatrix}
    x(0), \ldots, x(n_\D-2)
    \end{bmatrix}, \quad & X^+ = \begin{bmatrix}
        x(1), \ldots, x(n_\D-1)
    \end{bmatrix}, \\
U = \begin{bmatrix}
    u(0), \ldots, u(n_\D-2)
\end{bmatrix}, \quad & Y = \begin{bmatrix}
    y(0), \ldots, y(n_\D-2)
\end{bmatrix},
\end{aligned}
\]
and the lifted state can be computed as 
\[
\begin{aligned}
Z &= \begin{bmatrix}
    \Phi(x(0)), \ldots, \Phi(x(n_\D-2))
\end{bmatrix}, \\
Z^+ &= \begin{bmatrix}
    \Phi(x(1)), \ldots, \Phi(x(n_\D-1))
\end{bmatrix}.
\end{aligned}
\]
Finally, the matrices $A,B,C$ and $D$ can be obtained using least-squares approximations:
\begin{equation}
\label{eqn:EDMD}
    \begin{aligned}
        (A,B) &\in \argmin_{(A,B)} \|Z^+ -AZ -BU\|_F^2, \\
        (C,D) & \in \argmin_{(C,D)} \|Y-CZ-DU\|_F^2.
    \end{aligned}
\end{equation}

The accuracy of the parametric linear model \cref{eqn:Koop-linear} depends critically on the choice of the lifting functions $\Phi(\cdot)$, and an improper choice can lead to large bias errors~\cite{haseli2021learning}. Common function classes utilized to select $\Phi(\cdot)$ include polyharmonic splines, thin plate splines and Gaussian kernel~\cite{mauroy2020koopman}. However, a systematic approach for deciding the suitable parameters of lifting functions is highly non-trivial, while using deep neural networks (DNN) to learn lifting functions requires a large amount of data with a computationally expensive offline training process \cite{shi2022deep, han2020deep}. 

\subsection{Willems' fundamental lemma for CF-LCC system}
\label{subsec:Willems}
Willems' fundamental lemma is established for LTI systems, and thus the CF-LCC system \eqref{eqn:CF-LCC} needs to be linearized around an equilibrium state. From the behavioral (\emph{i.e.}, trajectory) perspective, the whole trajectory space of the linearized system can be represented as a linear combination of its rich enough offline trajectories.

We represent the equilibrium velocity for each vehicle as $v^*$ (all vehicles move in the same velocity) and the~equilibrium spacing as $s_i^*$ (which may vary from different vehicles). Then, we can define the velocity error and spacing error for each vehicle as $\tilde{v}_i(k) = v_i(k)-v^*, \tilde{s}_i(k) = s_i(k)-s_i^*$ and error state of the CF-LCC system can be derived as 
\[
\tilde{x}(k) = \begin{bmatrix}
    \tilde{s}_1(k), \tilde{v}_1(k), \ldots, \tilde{s}_n(k), \tilde{v}_n(k)
\end{bmatrix}^\tr \in \mathbb{R}^{2n}.
\]
The equilibrium velocity $v^*$ can be estimated from the past velocity trajectory of the head vehicle while the equilibrium spacing of HDVs is non-trivial to obtain \cite{wang2022implementation}. Thus, the measurable states are velocity errors of all vehicles and the spacing error of the CAV which can be designed. The output of the linearized system becomes 
\[
\tilde{y}(k) = \begin{bmatrix}
    \tilde{v}_1(k), \tilde{v}_2(k), \ldots, \tilde{v}_n(k), \tilde{s}_1(k)
\end{bmatrix}^\tr \in \mathbb{R}^{n+1}.
\]

Then, the dynamics of the linearized CF-LCC system in the error state space are  
\begin{equation}
\label{eqn:CF-LCC-linear}
\begin{aligned}
\tilde{x}(k+1)  = \tilde{A} \tilde{x}(k) + \tilde{B} u(k), \qquad 
\tilde{y}(k)  = \tilde{C} \tilde{x}(k),
\end{aligned}
\end{equation}
where the matrices $\tilde{A}, \tilde{B}$ and $\tilde{C}$ can be found in \cite{wang2023deep}.

We recall a persistent excitation condition for collecting rich enough offline trajectories.
\begin{definition}[Persistently exciting]
    The length-$T$ data sequence $\omega = \col(\omega(0), \ldots, \omega(T-1))$ is persistently exciting (PE) of order $L$ if its associated Hankel matrix 
    \[ 
\mathcal{H}_L(\omega) = \begin{bmatrix}
    \omega(0) & \omega(1) & \cdots &  \omega(T-L) \\
    \omega(1) & \omega(2) & \cdots & \omega(T-L+1) \\
    \vdots & \vdots & \ddots & \vdots \\
    \omega(L-1) & \omega(L) &\cdots  & \omega(T-1)
\end{bmatrix}
\]
has full row rank. 
\end{definition}

Then, with the pre-collected input-output data sequence $u_\D \!= \!\col(u(0), \ldots, u(n_\D-1)), y_\D \!=\! \col(y(0), \ldots, y(n_\D -1))$ of~\cref{eqn:CF-LCC-linear}, we have the following Willems' fundamental lemma.

\begin{lemma}[Willems' fundamental lemma \cite{willems2005note}]
\label{lemma:Willems}
    Consider the LTI system \cref{eqn:CF-LCC-linear}. Suppose the input trajectory $u_\D$ is persistently exciting of order $L+2n$. Then, an input-output data sequence $\col(u, y)\in \mathbb{R}^{(n+3)L}$ is a valid trajectory of \cref{eqn:CF-LCC-linear} if and only if there exists $g \in \mathbb{R}^{n_\D-L+1}$ such that 
    $
    \col(\mathcal{H}_L(u_\D), \mathcal{H}_L(y_\D)) g = \col(u, y).
    $
\end{lemma}

\vspace{2pt}

We can construct a data-driven representation for the linearized system \cref{eqn:CF-LCC-linear} using \Cref{lemma:Willems}. We denote $u_\ini = \col(u(k-T_\ini), \ldots, u(k-1))$ and $u_\f=\col(u(k), \ldots, u(k+N-1))$ as the most recent past length-$T_\ini$ input trajectory and the future length-$N$ input trajectory and $L = T_{\ini}+N$ (similarly for $y_\ini, y_\f$). We further divide the Hankel matrix into its first $T_\ini$ rows (\emph{i.e.}, $U_\p, Y_\p$) and the last $N$ rows (\emph{i.e.}, $U_\f, Y_\f$), which is 
\begin{equation}
\label{Hankel-partition}
\begin{bmatrix}
    U_\p \\ U_\f
\end{bmatrix}
:= \mathcal{H}_L(u_\D), \quad 
\begin{bmatrix}
Y_\p \\ Y_\f
\end{bmatrix}
:= \mathcal{H}_L(y_\D).
\end{equation}
Then, $\col(u_\ini, y_\ini, u_\f, y_\f)$ is a valid trajectory of \cref{eqn:CF-LCC-linear} if and only if there exists $g \in \mathbb{R}^{n_\D-T_\ini-N+1}$ which satisfies
\begin{equation}
\label{eqn:DD-CF-LCC-Linear}
\col(U_\p, Y_\p, U_\f, Y_\f)g = \col(u_\ini, y_\ini, u_\f, y_\f).
\end{equation}
We further note that $y_\f$ is unique for any $(u_\ini, y_\ini, u_\f)$ if $T_\ini \ge 2n$ and $\col(U_\p, Y_\p, U_\f, Y_\f)$ can be considered as a trajectory library for \eqref{eqn:CF-LCC-linear} (\emph{i.e.}, each column is a valid trajectory of the system).
\begin{remark}[Local linearization and the model update] It is clear that the data-driven representation \cref{eqn:DD-CF-LCC-Linear} depends on the linearized system \cref{eqn:CF-LCC-linear}. Thus, \cref{eqn:DD-CF-LCC-Linear} must be recomputed whenever the traffic equilibrium state changes. Furthermore, since the equilibrium spacing of HDVs is unknown, the measurable spacing information of HDVs is not used. \hfill $\square$
\end{remark}

\subsection{Problem statement}
In this work, we aim to first (approximately) represent the CF-LCC system \eqref{eqn:CF-LCC} with a dictionary-free data-driven representation for the Koopman linear model \cref{eqn:Koop-linear}. We then integrate the linear representation with the predictive control to formulate the proposed \method{DF-KMPC}, which is of the form 
\begin{subequations}
    \label{eqn: DD-K-pred}
    \begin{align}
        \min_{g, u \in \mathcal{U}, y \in \mathcal{Y}} \quad & \|y- y_\R \|_Q + \|u\|_R \\
        \mathrm{subject~to} \quad & \Bar{H}^* g = \col(u_\ini, y_\ini, u, y),  \label{eqn:DD-K}
    \end{align}
\end{subequations}
where $\bar{H}^*$ is similar to a Hankel matrix encoding a dictionary-free Koopman linear representation, $\mathcal{U}, \mathcal{Y}$ are the input, output constraints. The reference trajectory are represented as $y_r := \mathbb{1} \otimes \col(s_\R, v_\R) $ where $\mathbb{1} \! \in \! \mathbb{R}^{Nn}$ is a column vector with all elements equal to~$1$ and $s_r, v_r$ denote the reference spacing and velocity respectively.

The dictionary-free Koopman linear representation \cref{eqn:DD-K} is built on our recent work \cite{shang2024willems}. Indeed, the theoretical work~\cite{shang2024willems} reveals that one can directly construct a data-driven representation for the Koopman linear model via Willems' fundamental lemma when an exact Koopman embedding exists. However, whether the mixed traffic system has an exact Koopman linear embedding is still an open problem. When the Koopman linear model is inexact, its encoded input and output trajectories are different from the true behaviors of the CF-LCC system,  and thus the results in~\cite{shang2024willems} can not be directly applied. We tackle this challenge in two steps. We first formulate the construction of the optimal Koopman linear model as an optimization problem from the behavioral perspective. We further~relax this problem and modify the constraints based on system identification technologies to obtain the data-driven representation. We~then adapt an iterative algorithm \cite{shang2024convex} to approximately compute the Hankel matrix $\Bar{H}^*$ using for online predictive control.  

\begin{remark}[Deep Hankel matrix and no explicit lifting]
    We construct the data-driven representation \cref{eqn:DD-K} for the Koopman linear model \cref{eqn:Koop-linear} via processing the input-output trajectories of the CF-LCC system \cref{eqn:CF-LCC} and do not require explicit lifted states. One key insight from \cite{shang2024willems} is that the Hankel matrix  $\bar{H}^*$ should be sufficiently deep.  The resulting  $\bar{H}^*$ can be viewed as an approximated Koopman linear model around the operating region (\emph{i.e.}, the region where data is collected). Thus, updating the equilibrium state of the mixed traffic system as in \cite{wang2023deep, shang2024decentralized} is not required in~\cref{eqn: DD-K-pred}. 
\end{remark}

\section{Predictive Control with Data-driven Koopman Linear Models}
\label{sec:data-driven-Koopman}
We here present a dictionary-free Koopman linear model for the CF-LCC system in two cases: 1)~an~exact~Koopman linear model exists and 2) the Koopman linear model is inexact. We then introduce \method{DF-KMPC} for CAV~control.

\subsection{CF-LCC system with exact Koopman linear embedding} 
We here illustrate the construction of data-driven representation for an exact Koopman linear embedding for the CF-LCC system via adapting the results in \cite{shang2024willems}. We first introduce a notion of \textit{lifted excitation} for the pre-collected input-state-output trajectories $u_\D, x_\D$ and $y_\D$. 
\begin{definition}[Lifted excitation]
    Suppose there exists an exact Koopman linear embedding \cref{eqn:Koop-linear} for the CF-LCC system~\cref{eqn:CF-LCC}. We say the input-output trajectory $u_\D, y_\D$ collected from \cref{eqn:CF-LCC} provides lifted excitation of order $L$, if 
    \[
    H_K := \begin{bmatrix}
        u(0) & u(1) & \cdots & u(n_\D-L) \\
        u(1) & u(2) & \cdots & u(n_\D-L+1) \\
        \vdots & \vdots & \ddots & \vdots \\
        u(L-1) & u(L) & \cdots & u(n_\D-1) \\
        \Phi(x(0)) & \Phi(x(1)) & \cdots & \Phi(x(n_\D-L))
    \end{bmatrix}
    \]
    has full row rank.
\end{definition}

This definition provides a sufficient condition for the collected trajectory from \cref{eqn:CF-LCC} to be rich enough to formulate the trajectory space of the Koopman linear embedding~\cref{eqn:Koop-linear}. We then construct a data-driven representation for the CF-LCC system (equivalent to the exact Koopman linear model~\cref{eqn:Koop-linear}). We divide the Hankel matrix for the input and output trajectory the same as in \cref{Hankel-partition} and obtain $U_\p, U_\f, Y_\p, Y_\f$. 
\begin{theorem}
    Suppose there exists an exact Koopman linear embedding \cref{eqn:Koop-linear} for the CF-LCC system \cref{eqn:CF-LCC}. We collect the input-output trajectory $u_\D, y_\D$ that has lifted excitation of order $L$, where $L = T_\ini +N$ and $T_\ini \ge n_\z$. At time $k$, we denote the most recent length-$T_\ini$ input-output trajectory $\col(u_\ini, y_\ini)$ from \cref{eqn:CF-LCC} as 
    \[
    \begin{aligned}
    u_\ini & = \col(u(k-T_\ini), \ldots, u(k-1)), \\
    y_\ini & = \col(y(k-T_\ini), \ldots, y(k-1)).
    \end{aligned}
    \]
    For any future input $u_\f = \col(u(k), \ldots, u(k+N-1))$, the sequence $\col(u_\ini, y_\ini, u_\f, y_\f)$ is a valid length-$L$ trajectory of \cref{eqn:CF-LCC} if and only if there exists $g \in \mathbb{R}^{n_\D-T_\ini-N+1}$ such that 
    \[
    \col(U_\p, Y_\p, U_\f, Y_\f) g = \col(u_\ini, y_\ini, u_\f, y_\f).
    \]
\end{theorem}

\vspace{3pt}

The proof can be adapted from \cite[Theoerem 3]{shang2024willems}.  This result provides a \textit{dictionary-free} data-driven representation for the nonlinear CF-LCC system \cref{eqn:CF-LCC} by considering the trajectory space of the Koopman linear embedding. An important feature is that no explicit lifting is required. Instead, the online initial trajectory should be long enough, depending on the hidden dimension of the Koopman linear model (throughout the paper, we consider $T_\ini = n_\z$). Thus, the Hankel matrix should have sufficient depth. 

\subsection{CF-LCC system with inexact Koopman linear model}
As the CF-LCC system may not have an exact Koopman linear embedding, we here construct a data-driven representation for an approximated Koopman linear model of the nonlinear system. Different from the EDMD approach~\cref{eqn:EDMD} using pre-selected lifting functions, we seek the optimal Koopman linear model from a behavioral perspective, which bypasses the process of selecting lifting functions. 

We consider the trajectory space of the Koopman linear model and aim to minimize the distance between the collected data $u_\D, y_\D$ of \cref{eqn:CF-LCC} and the trajectory space. The optimization problem can be formulated as follows
\begin{equation}
    \label{eqn:opt-koop-model}
    \begin{aligned}
        \min_{\substack{A, B, C ,D, \\ \bar{y},\Phi \subset \mathcal{F}}} \quad & \| \Bar{y} -  y_\D \|_2 \\
        \mathrm{subject~to} \quad &  \col(u_\D, \Bar{y}) \in \mathcal{B}(A,B,C,D, \Phi),   
    \end{aligned}
\end{equation}
where $\mathcal{F}$ is a given function class and $\mathcal{B}_L$ denotes the length-$L$ trajectory space of the Koopman linear model 
\[
\mathcal{B}_L = \left\{\!\begin{bmatrix}
            u\\y
        \end{bmatrix}  \mid   \exists x_0 \! \in \! \mathbb{R}^{2n}, \; \cref{eqn:Koop-linear} \; \text{holds with} \; z(0)\! = \! \Phi(x_0) \!\right\}\!.
\]
We note that, as the input of the system is accurate, we focus on optimizing the output trajectory which represents the response of the dynamic system.

We further reformulate \eqref{eqn:opt-koop-model} by organizing input-output trajectory in a Hankel matrix which can be considered as a trajectory library. We also replace the constraint of the trajectory space of the Koopman linear model with a set of constraints arising from linear system identification techniques. Motivated by our work \cite{shang2024convex}, we relax \cref{eqn:opt-koop-model} as
\begin{subequations} \label{eqn:opt-Koop-DD}
    \begin{align} 
\min_{\bar{H}_y, K} \quad & \|H_y - \bar{H}_y\|_F \nonumber\\
\mathrm{subject~to} \quad & \textrm{rank}(\bar{H}) = mL+n_\z  \label{eq:SVD-Dom-1},\\
& \bar{Y}_\f = K \ \col(U_\p, \bar{Y}_\p, U_\f), \label{eq:SVD-Dom-3}\\
& K = \begin{bmatrix}
    K_p & K_f
\end{bmatrix}, \ K_f \in \mathcal{L} \label{eq:SVD-Dom-4}, \\ 
& \bar{H}_y \in \mathcal{H}, \label{eq:SVD-Dom-2}
\end{align}
\end{subequations}
where $H_y, \bar{H}_y$ and $\bar{H}$ represent $ \col(Y_\p, Y_\f)$, $ \col( \bar{Y}_\p, \bar{Y}_\f)$ and $\col(U_\p, \bar{Y}_\p, U_\f, \bar{Y}_\f)$, respectively, (see the partition in \cref{Hankel-partition}) and $\mathcal{L}$ represents the lower-block triangular matrix. The key insight is that the input-output trajectory has some specific properties if it comes from an LTI system. The constraint~\cref{eq:SVD-Dom-1} represents a low-rank constraint and \cref{eq:SVD-Dom-2} implies $\bar{H}_y$ needs to satisfy the Hankel structure \cite{markovsky2008structured}. Constraints \cref{eq:SVD-Dom-3} indicates that the future output $\bar{Y}_\f$ is a linear combination of past data $\col(U_\p, \bar{Y}_\p)$ and the future input $U_\f$~\cite{favoreel1999spc}, while \cref{eq:SVD-Dom-4} enforces the causality by requiring $K_f$ to be a lower-block triangular matrix \cite{sader2023causality}. The optimization problem \eqref{eqn:opt-Koop-DD} is a relaxation of \eqref{eqn:opt-koop-model} as all feasible solution $\bar{y}$ in \eqref{eqn:opt-koop-model} is feasible for \eqref{eqn:opt-Koop-DD} after reformulation. 

The optimization problem \cref{eqn:opt-Koop-DD} is still challenging to solve with all constraints. We address them in an alternating minimization process by adapting an iterative algorithm from \cite{shang2024convex}. We first consider the constraint \cref{eq:SVD-Dom-1} and relax \eqref{eqn:opt-Koop-DD} as 
\begin{equation}
\label{eqn:low-rank-approx}
    \begin{aligned} 
\min_{\bar{H}_y} \quad & \|H_y - \bar{H}_y\|_F \\
\mathrm{subject~to} \quad & \textrm{rank}(\bar{H}) = mL+n_\z.
\end{aligned}
\end{equation}
An analytical solution of  \cref{eqn:low-rank-approx} can be provided via singular value decomposition and we denote the mapping from $H_y$ (parameter of the optimization problem \cref{eqn:low-rank-approx}) to its optimal solution $H_{y_1}$ as $\Pi_\Lo$. We then take the obtained low-rank approximation $H_{y_1}$ as the parameter of the optimization
\begin{equation}
\label{eqn:row-causal}
    \begin{aligned} 
\min_{\bar{H}_y, K} \quad & \|H_{y_1} - \bar{H}_y\|_F \\
\mathrm{subject~to} \quad & \bar{Y}_\f = K \ \col(U_\p, Y_{\p_1}, U_\f), \\
& K = \begin{bmatrix}
    K_p & K_f
\end{bmatrix}, \ K_f \in \mathcal{L},
\end{aligned}
\end{equation}
which tackles constraints \cref{eq:SVD-Dom-3}, \cref{eq:SVD-Dom-4} and $\col(Y_{\p_1}, Y_{\f_1})\! :=\! H_{y_1}$. The problem~\cref{eqn:row-causal} also has an analytical solution \cite{sader2023causality} and we denote the mapping from the parameter $H_{y_1}$ to the optimal solution $H_{y_2}$ as $\Pi_\C$. Finally, we project  $H_{y_2}$ to a Hankel matrix set via averaging skew-diagonal elements \cite{yin2021low} and represent the projector and the resulting Hankel matrix as $\Pi_\h$ and $H_{y_3}$. The above steps are repeated iteratively until the result converges practically. We summarize the overall procedure in \Cref{alg:iter-SLRA}. The data-driven representation of an approximated Koopman linear model becomes $\bar{H}^* = \col(U_\p, Y_\p^*, U_\f, Y_\f^*)$ where $\col(Y_\p^*, Y_\f^*) := H_y^*$.
\begin{algorithm}[t]
\caption{Iterative Hankel-Koopman Construction}
    \label{alg:iter-SLRA}
  \begin{algorithmic}[1]
    \Require 
    $U_\p$,  $U_\f$, $Y_\p$, $Y_\f$, $n_\z$, $\epsilon$
    \State $H_{y} \gets \col(Y_\p, Y_\f)$, $H_{y_3} \gets H_{y}$;
    \Repeat
    \State $H_{y_1} \gets \Pi_\Lo(H_{y_3})$ \quad (\texttt{Low-rank approx});
    \State $H_{y_2} \!\gets\!\! \Pi_\C (H_{y_1})$ \quad (\texttt{Causality~proj});
    \State $H_{y_3} \!\gets\!\! \Pi_\h (H_{y_2})$ \quad (\texttt{Hankel~proj});
    \Until{$\|H_{y_1} -H_{y_3}\|_F \le \epsilon \|H_{y_1}\|_F$}
    \Ensure $H_y^* = H_{y_1}$
\end{algorithmic}
\end{algorithm}

\begin{algorithm}[t] 
	\caption{\method{DF-KMPC}}
	\label{Alg:DK-MPC}
	\begin{algorithmic}[1]
		\Require
		Pre-collected offline data $(u_{\textrm{d}}, y_{\textrm{d}})$, initial time step $k_0$, terminal time step $k_f$;
		\State Construct data Hankel matrices for input, and output as $U_{\textrm{P}}, U_{\textrm{F}}, Y_{\textrm{P}}, Y_{\textrm{F}}$;
            \State Develop approximated data-driven representation $\bar{H}^*$ for optimal Koopman linear embedding via \Cref{alg:iter-SLRA};
		\State Initialize the most recent past traffic data $(u_{\textrm{ini}},y_{\textrm{ini}})$ before the initial time $k_0$;
		\While{$k_0 \leq k \leq k_f$}
		\State Solve~\eqref{eqn:DD-K-pred-proj} for optimal future control input sequence \textcolor{white}{\, \, \, }$u^*=\textrm{col}(u^*(k),u^*(k+1),\ldots,u^*(k+N-1))$;
		\State Apply the input $u(k) \leftarrow u^*(k)$ to the CAV;
		\State $k \leftarrow k+1$;
            \State Update past traffic data $(u_{ \textrm{ini}},y_{\textrm{ini}})$;
		\EndWhile
	\end{algorithmic}
\end{algorithm}

\subsection{Dictionary-free Koopman model predictive control}
After obtaining the dictionary-free Koopman linear representation, the problem \eqref{eqn: DD-K-pred} is well defined. However, as the exact Koopman linear representation may not exist, it is possible that the trajectory of the nonlinear system is not included in the trajectory space of the data-driven representation. That may lead to an infeasible solution. Thus, we project the most recent input-output trajectory $\col(u_\ini, y_\ini)$ to the range space of $\col(U_\p, Y_\p^*)$ which ensures the feasibility of the initial condition. We denote the projector as $\Pi_\ini$ and the final form of \eqref{eqn: DD-K-pred} becomes
\begin{subequations}
\label{eqn:DD-K-pred-proj}
    \begin{align}
        \min_{g, u \in \mathcal{U}, y \in \mathcal{Y}} \quad & \|y- y_\R \|_Q + \|u\|_R \\
        \mathrm{subject~to} \quad & \Bar{H}^* g = \col(\Pi_\ini(\col(u_\ini, y_\ini)), u, y).
    \end{align}
\end{subequations}

We note that most of existing data-driven control approaches (\emph{e.g.}, \cite{wang2023deep, shang2024decentralized}) for the mixed traffic system require 1) the selection of regularizer on $g$ to implicitly identify the system model \cite{dorfler2022bridging}, and 2) the extra slack variable $\sigma_y$ adding on $y_\ini$ to ensure the feasibility of \cref{eqn:DD-K-pred-proj}. Our approach does not require either of them as obtaining the data-driven Koopman linear model is already an implicit system identification and we project $\col(u_\ini, y_\ini)$ to the trajectory space $\mathcal{B}_{T_{\ini}}$. The overall procedure of the \method{DF-KMPC} is detailed in \Cref{Alg:DK-MPC}.

\section{Numerical Experiments}
\label{sec:results}
In this section, we conduct nonlinear traffic simulations to compare the performance of predictive control with different linear models: 1) our proposed \method{DF-KMPC} with dictionary-free Koopman representation (\method{DF-K}), 2) predictive control with the standard Koopman linear model \cref{eqn:Koop-linear} obtained from EDMD \cref{eqn:EDMD} (\method{EDMD-K}) and 3) predictive control with Deep NN Koopman representation (\method{DNN-K}) in \cite{shi2022deep}.

\subsection{Experiment setup}
\begin{figure}[t]
    \centering
    \includegraphics[width=0.3\textwidth]{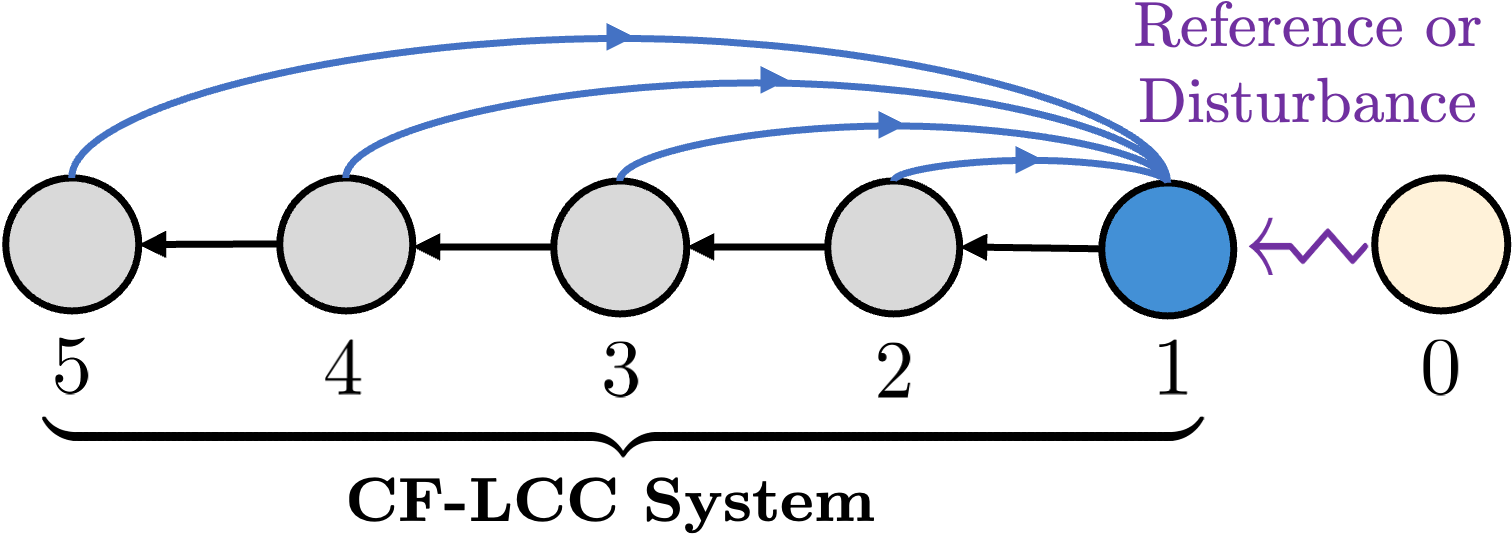}
    \caption{Simulation scenario. The CF-LCC system consists of $5$ vehicles, where the yellow node, blue node, and grey nodes represent the head vehicle, the CAV, and other HDVs, respectively.}
    \label{fig:CF-LCC-Sim}
\end{figure}

Similar to \cite{wang2023deep,shang2024decentralized}, we consider a CF-LCC system with one CAV and four HDVs; see \Cref{fig:CF-LCC-Sim}. The velocity of the head vehicle is taken either as a reference or a disturbance for the CF-LCC system. We model the car-following behaviors of HDVs by the nonlinear OVM model proposed in \cite{bando1995dynamical}.

We use the following parameters for constructing the data-driven linear models and synthesizing controllers:
\begin{enumerate}
    \item \textit{Offline data collection:} we set the time step $\Delta t=0.05\,\mathrm{s}$. The initial state of each vehicle is randomly selected with a uniform distribution of $v_i \in [10, 20]$ and $s_i \in [15, 25]$ for $i = 1,\ldots 5$. The $u_1, v_0$ are generated by uniform distributed signals of $\mathbb{U}[-5, 5]$ and $\mathbb{U}[10, 20]$ respectively. We collect a single trajectory with length-1200 for \method{DF-K} and choose $T_\ini$ to be $40$ to implicitly represent a $40$-dimension Koopman linear model. To obtain the Koopman linear model~with EDMD, we use 100 trajectories of length 1200. The lifting functions are chosen to be the state of \cref{eqn:CF-LCC} and $30$ thin plate spline radial basis functions with the form $\phi(x) = \|x-x_0\|_2^2\log(\|x-x_0\|_2)$, whose center $x_0$ is random sampled from a uniform distribution $[5, 15]^5 \times [10, 20]^5$. To learn lifting functions with~DNN, we simulate 2000 trajectories with 500 steps and utilize 2 hidden layers with 32 units to learn 40 lifting functions. Thus, \method{DF-K}, \method{EDMD-K} and \method{DNN-K} all construct a 40-dimension (implicit) Koopman linear model with  $1.2\times 10^3, 1.2\times 10^5$ and $10^6$ data points,  respectively. 
    \item \textit{Online predictive control:} we set the input and output constraints for the acceleration and spacing of the CAV as $\mathcal{U} := \{u \! \in \! \mathbb{R}^{2N} \mid a_{\min} \! \le \! u_i \! \le  \! a_{\max}, \, \text{for}$ $\text{$i = 1, 3,\ldots, 2N-1$}\}$ and $\mathcal{Y} := \{y \in \mathbb{R}^{10N} \mid s_{\min} \le s_i \le s_{\max} \, \text{for $i = 1, 11,\ldots, 10N-9$}\}$. The prediction horizon is set to $N=50$ and we have $R = 0.1I_N$, $Q = I_{Nn} \otimes \mathrm{diag}(Q_v, w_s)$ where $Q_v, w_s$ are varying for scenarios with different control objectives. The limitation for the acceleration of CAV is set as $a_{\max} = 2\,\mathrm{m/s^2}$ and $a_{\min} = -5\,\mathrm{m/s^2}$ and its spacing constraint are $s_{\max} = 40\,\mathrm{m}$ and $s_{\min} = 5 \, \mathrm{m}$.
\end{enumerate}

\begin{figure}[t]
\centering
\setlength{\abovecaptionskip}{2pt}
\hspace{-2mm}\subfigure[HDV]{\includegraphics[width=0.24\textwidth]{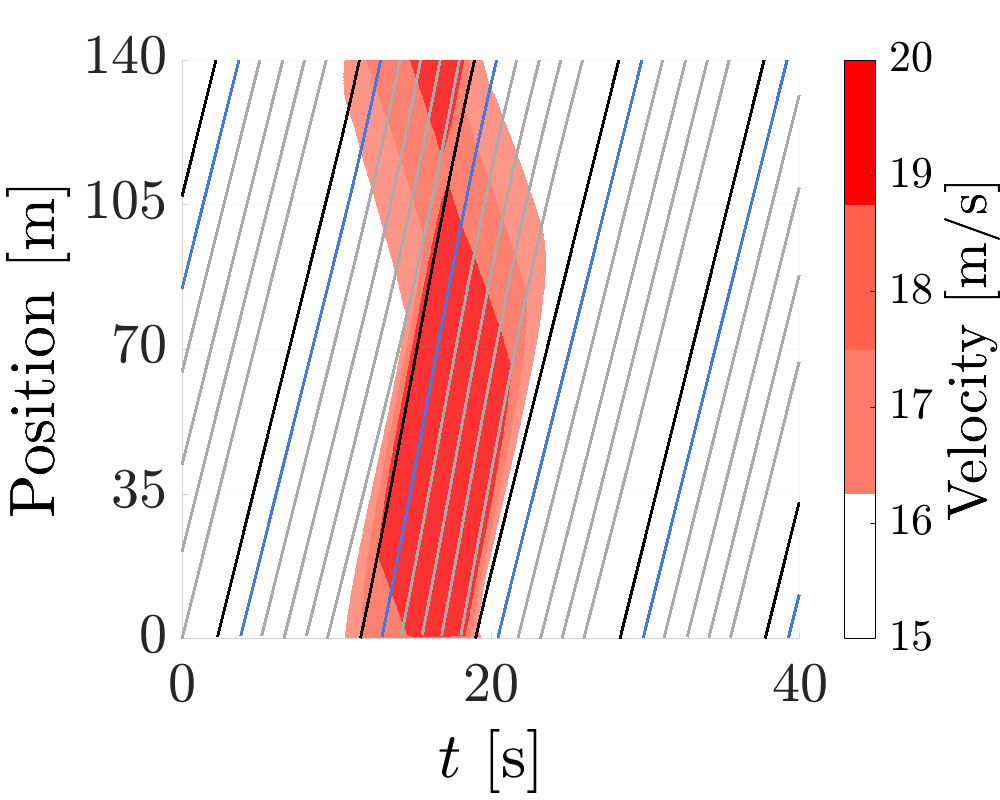}\label{subfig:HDV}}
\subfigure[\method{DF-K}]{\includegraphics[width=0.24\textwidth]{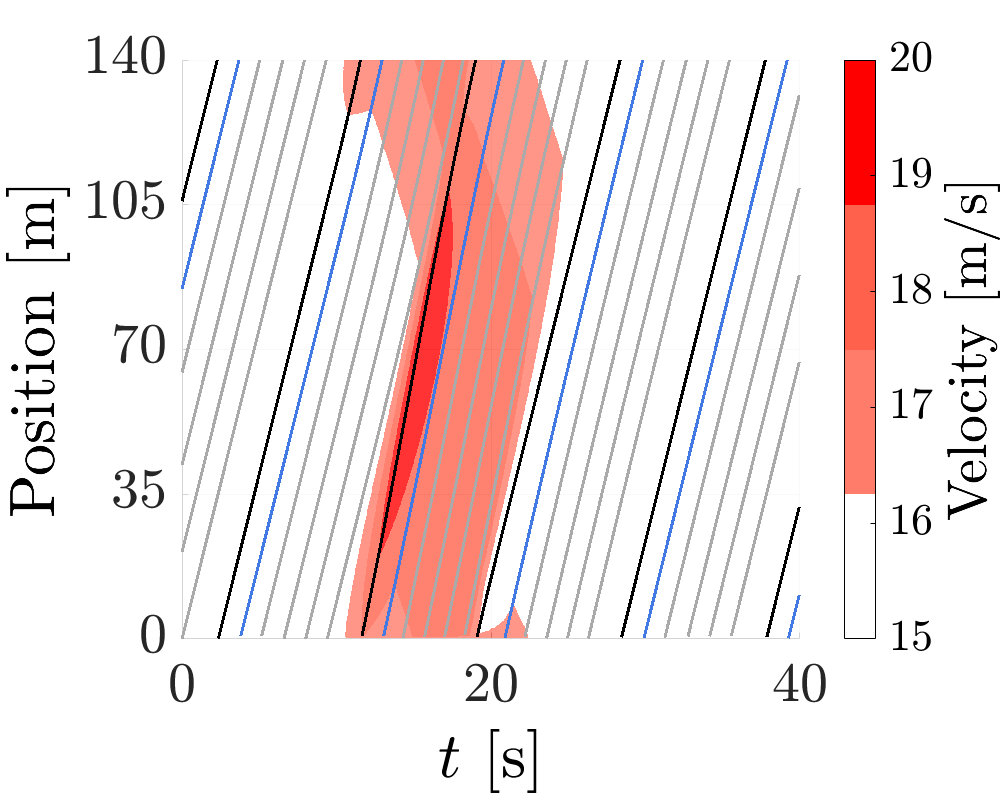} \label{subfig:DDK}} \\ 
\hspace{-4mm} \subfigure[\method{EDMD-K}]{\includegraphics[width=0.24\textwidth]{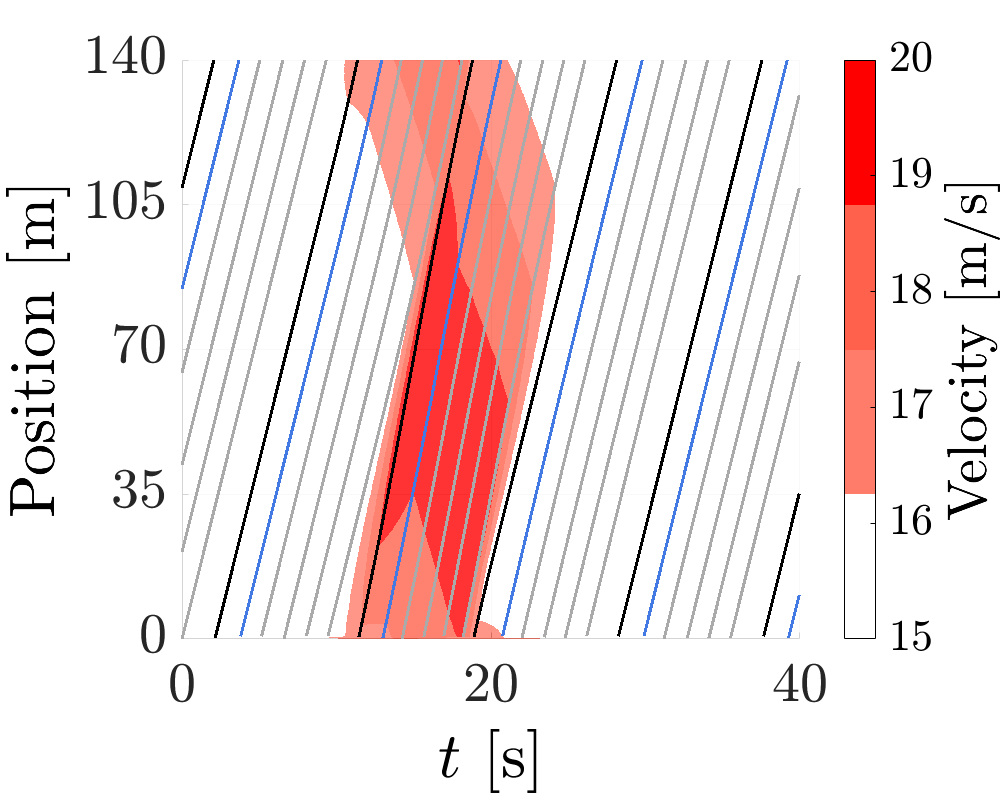}\label{subfig:EDMDK}}
\subfigure[\method{DNN-K}]{\includegraphics[width=0.24\textwidth]{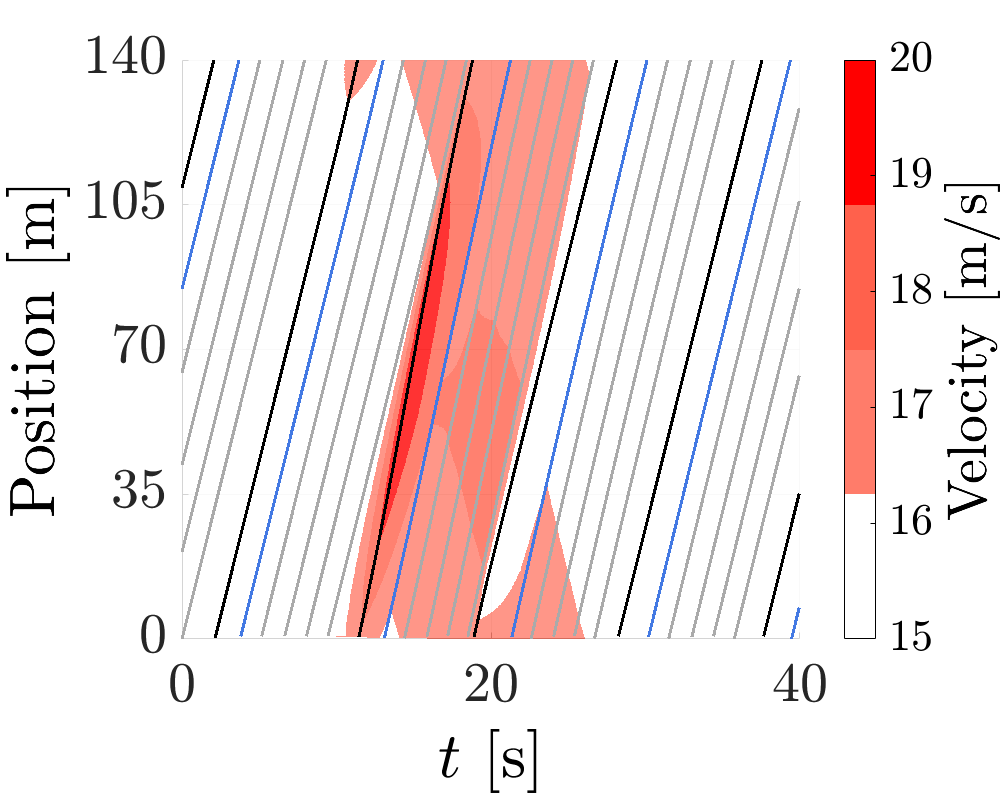}\label{subfig:DNN-K}}
\caption{Schematic of trajectories and velocities of vehicles in Experiment~A. The black line and blue line represent the position of the head vehicle and the CAV, respectively. The darker the color, the higher the velocity. The disturbance applied on the head vehicle starts at $t = 10 \, \mathrm{s}$ and ends at $t = 20 \, \mathrm{s}$. (a) All vehicles are HDVs. (b)-(d) The CAV utilizes predictive controllers with different linear representations. }
\label{fig:time-space}
\vspace{-5mm}
\end{figure}

\subsection{Traffic wave mitigating and trajectory tracking}
\textbf{Experiment A:} We first compare the performance of mitigating traffic waves by utilizing different linear representations in predictive control. We consider a ring road scenario with circumference $140\,\mathrm{m}$. The head vehicle is under the perturbation of a sine wave that is $v_0(k) = 15 + 5\sin(\pi k/200)$ and the predictive controller is required to regulate the traffic wave while maintaining the spacing between vehicles (\emph{i.e.}, $s_\R = 20\,\mathrm{m}$, $v_\R = 15\,\mathrm{m/s}$). 

The results are displayed in \Cref{fig:time-space}. From $t = 10\,\mathrm{s}$ to $t = 20\,\mathrm{s}$, a traffic wave is gradually generated due to the perturbation of the head vehicle (see the black line and its associated dark red area). It is obvious that when all vehicles follow the human driving behaviors, the traffic wave propagates along the vehicle chain without vanishing (see the dark red area in a parallelogram shape in \Cref{subfig:HDV}). When the CAV is equipped with the predictive controller utilizing \method{DF-K} or \method{DNN-K}, it can effectively prevent and dampen the propagation of the traffic wave (see the dark red region disappear after the blue line in \Cref{subfig:DDK} and  \Cref{subfig:DNN-K}). Although the predictive controller using \method{EDMD-K} can also mitigate the traffic wave as the dark red region in \Cref{subfig:EDMDK} is smaller than \Cref{subfig:HDV}, it is much larger than that in \Cref{subfig:DDK} and  \Cref{subfig:DNN-K}. That illustrates a random or improper choice of lifting functions can lead to an inaccurate approximation, which deteriorates the control performance. On the other hand, the proposed \method{DF-K} approximates a data-driven representation for the optimal Koopman linear model, which improves the model accuracy and shows satisfied traffic wave mitigation performance.

\textbf{Experiment B:} We further demonstrate the trajectory tracking performance of predictive controllers with different linear models. We consider a highway scenario and utilize the real vehicle trajectory from the Next Generation SIMulation (NGSIM) program \cite{NGSIM}. We extract four vehicle velocity trajectories from the US-101 dataset which collects the traffic data on a freeway segment of US-101 and apply it as the velocity profile for the head vehicle. The predictive controller is then designed to track the velocity of the head vehicle while maintaining a safe distance between the CAV and the head vehicle (\emph{i.e.}, $v_\R$ is estimated from the past velocity trajectory of the head vehicle and $Q_v = 1, w_s = 0$). 

\begin{figure}[t]
\centering
\setlength{\abovecaptionskip}{2pt}
\subfigure[Velocity trajectorties]{\includegraphics[width=0.24\textwidth]{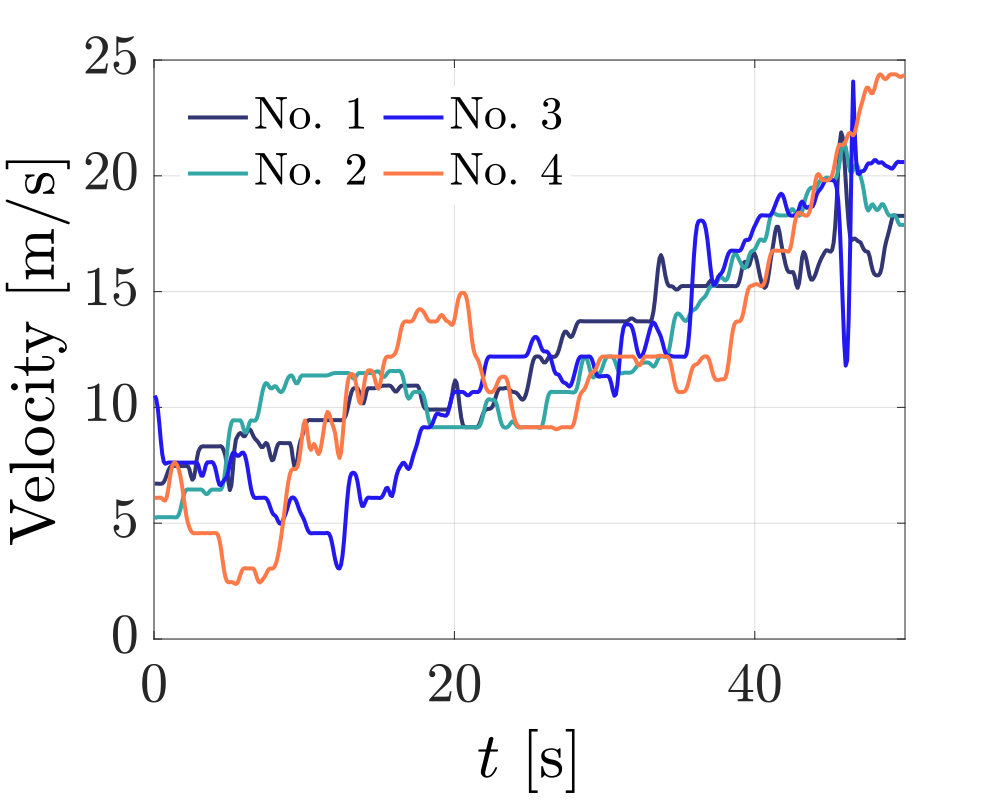}\label{subfig:NGSIM-traj}}
\hspace{-2mm}
\subfigure[Realization cost]{\includegraphics[width=0.24\textwidth]{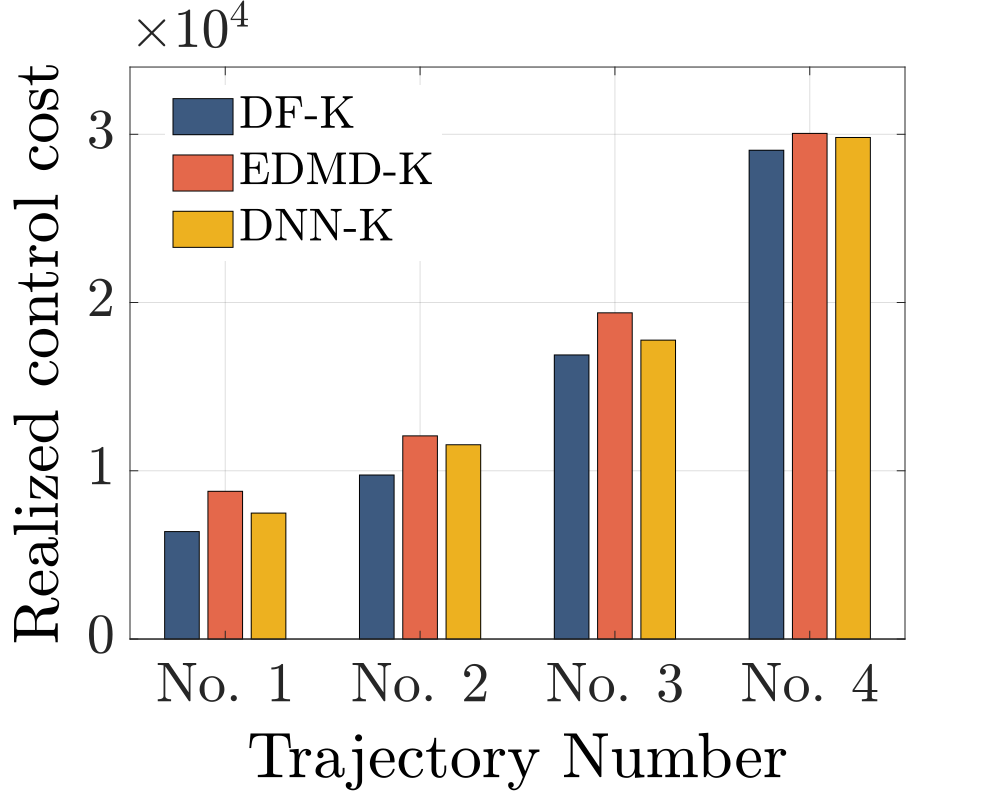}\label{subfig:real-cost}} 
\caption{Velocity profiles of head vehicle and realization cost of predictive controller with different linear models. (a) Real traffic velocity profiles of vehicles. (b) Realized control cost.}
\label{fig: traj-cost}
\vspace{-4mm}
\end{figure}

The four velocity profiles of the head vehicle are shown in \Cref{subfig:NGSIM-traj}, and realized control costs for different linear models are presented in \Cref{subfig:real-cost}. We consider the realized control cost that is $\|y^*-y_\R\|_Q^2+\|u^*\|_R^2$ where $u^*$ is the optimal control computed from the predictive controller and $y^*$ is the resulted actual trajectory. The realized control cost in \Cref{subfig:real-cost} illustrates $\method{EDMD-K} > \method{DNN-K} > \method{DF-K}$. We can clearly observe that the proposed \method{DF-K} outperforms the \method{EDMD-K} with a random choice of lifting functions. Although the \method{DNN-K} can provide comparable performance for some trajectories (see trajectories 3 and 4 in \Cref{subfig:real-cost}), it requires an order of magnitude more data ($1200$ \emph{v.s.} $10^6$) with much longer offline computational time (\emph{$15\,\mathrm{s}$ v.s. $2.1\times 10^4 \,\mathrm{s}$}).

\begin{figure}[t]
\centering
\setlength{\abovecaptionskip}{2pt}
\subfigure[\method{DF-K}]{\includegraphics[width=0.24\textwidth]{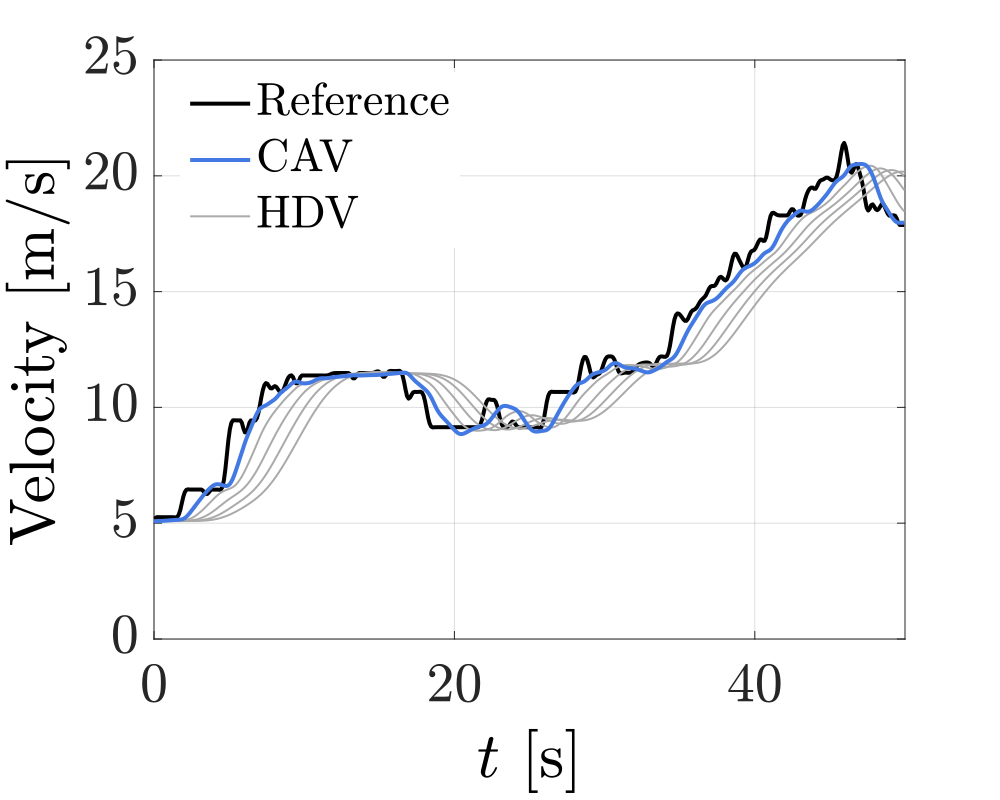}\label{subfig:track-DDK}}
\hspace{-2mm}
\subfigure[\method{EDMD-K}]{\includegraphics[width=0.24\textwidth]{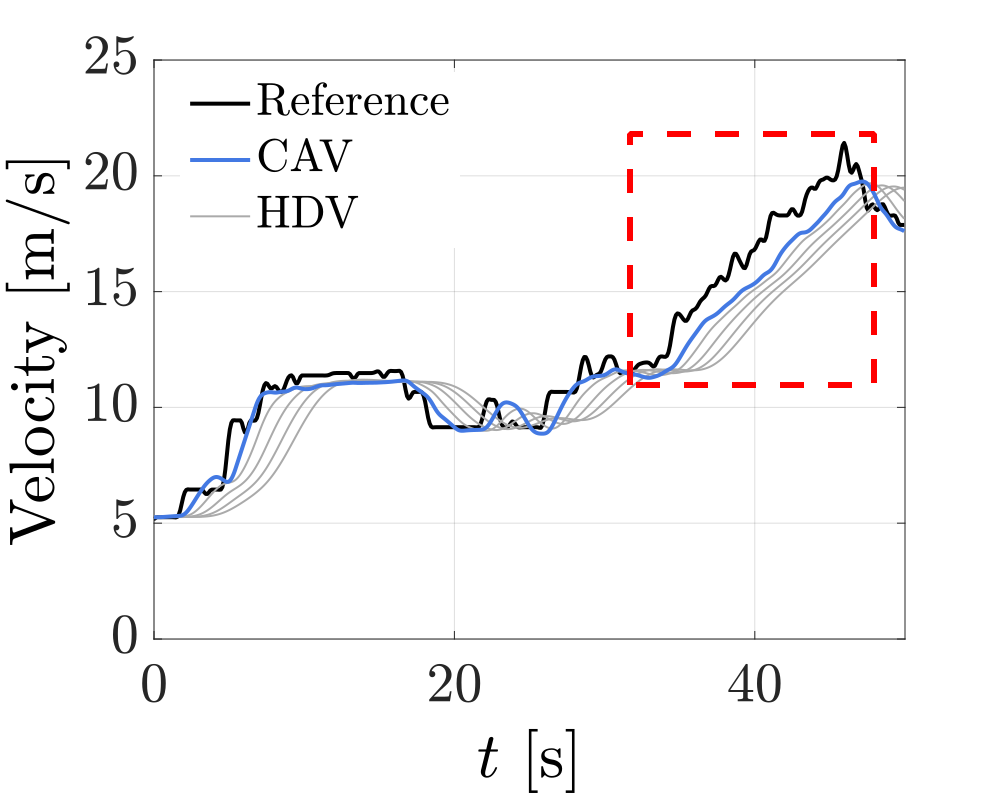}\label{subfig:track-EDMD}} 
\caption{Velocity profiles in Experiment B. The black profile denotes the head vehicle, while the blue profile and gray profiles represent the CAV and HDVs respectively. (a) The CAV utilizes the proposed \method{DF-KMPC} with \method{DF-K}. (b) The CAV utilizes a predictive controller with \method{EDMD-K}.}
\label{fig:tracking-perform}
\vspace{-2mm}
\end{figure}

\vspace{-2mm}
\Cref{fig:tracking-perform} presents the tracking performance of \method{DF-K} and \method{EDMD-K}. As shown in \Cref{subfig:track-DDK}, the CAV with the proposed \method{DF-K} can track the reference trajectory closely (see the blue and black curves in \Cref{subfig:track-DDK}). The tracking performance of CAV with \method{EDMD-K} can provide a satisfied tracking performance (\emph{i.e.}, follow the velocity of the head vehicle) in most of the time. However, a large deviation can occur (see the region inside the red dashed box) due to the relatively large approximation error induced by an improper set of lifting functions. The proposed \method{DF-K} approximates the optimal Koopman representation from the behavior perspective, bypassing the selection of the lifting functions

\begin{remark}
    For all data-driven approaches we compared in this paper (\emph{i.e.}, \method{DF-K}, \method{EDMD-K}, and \method{DNN-K}), their performance depends on the data collection region and the model can become relatively inaccurate when the system is out of the operating region. We note that the performance of \method{EDMD-K} and \method{DNN-K} may be improved with further tuning. However, there is no systemic procedure for the tuning process while we provide an efficient iterative algorithm to obtain the approximated data-driven representation for the optimal Koopman linear model.
\end{remark}

\section{Conclusion}
In this work, we have developed the \method{DF-KMPC} for CAV control in mixed traffic. We provided a systematic procedure to obtain an approximated data-driven representation for the optimal Koopman linear model. This data-driven representation eliminates the selection process of lifting functions and does not need to be updated when the traffic equilibrium state changes. Simulations with nonlinear traffic systems have validated the performance of the \method{DF-KMPC} in mitigating the traffic wave and tracking the desired velocity. Interesting future directions include considering the set of potential future velocities of the head vehicle, incorporating the spacing constraint of HDVs and testing the proposed dictionary-free \method{KMPC} in large-scale traffic simulation.
\label{sec:conclusion}

\bibliographystyle{IEEEtran}
\bibliography{ref}

\end{document}